\input epsf

\documentclass[10pt]{amsart}

\usepackage{amscd, amsmath, amsthm}

\def\mylabel#1{\label{#1}   \rlap{\hskip1cm\leftline{#1}}   }
\def\mylabel#1{\label{#1}   \proplabeL{#1} \hskip-3pt  }

\def\mylabel#1{\label{#1}}

\def\pd{\partial}

\newtheorem{theorem}{Theorem}[section]

\newtheorem{proposition}[theorem]{Proposition}

\theoremstyle{definition}     

\newtheorem{remark}[theorem]{Remark}

\numberwithin{equation}{section}

\begin{document}

\title[ On Stokes matrices of CY hypersurfaces]{ 
{On Stokes Matrices of Calabi-Yau hypersurfaces}
}

\author[]{}
\address{}
\email{}

\author[C.F. Doran]{ Charles F. Doran}
\author[S. Hosono]{Shinobu Hosono}

\address{
Department of Mathematics, 
University of Washington, Seattle, 
Washington 98195, U.S.A.}
\email{doran@math.washington.edu}

\address{
Graduate School of Mathematical Sciences, 
University of Tokyo, Komaba Meguro-ku, 
Tokyo 153-8914, Japan
}
\email{hosono@ms.u-tokyo.ac.jp}


\begin{abstract}
We consider Laplace transforms of the Picard-Fuchs differential 
equations of Calabi-Yau hypersurfaces and calculate their 
Stokes matrices. We also introduce two different types of 
Laplace transforms of Gel'fand-Kapranov-Zelevinski 
hypergeometric systems.
\end{abstract}

\maketitle


\section{ {\bf Introduction -- Laplace transforms of period integrals }}

Oscillatory integrals play important roles 
in Landau-Ginzburg theory of an N=2 supersymmetric 
field theory\cite{CV1,Du1}.  
They are characterized by certain differential 
equations with an irregular singulartity at infinity, and the Stokes 
matrices at infinity were key ingredients in the classification 
program of two dimensional N=2 supersymmetric field theories \cite{CV2}
and more recently in the homological mirror symmetry of Fano surfaces
\cite{HIV,AKO}.  

In this note, we consider the Landau-Ginzburg theories of 
Calabi-Yau hypersurfaces in weighted projective spaces 
and determine the Stokes matrices of the relevant 
oscillatory integrals. Because of a certain degeneration of the 
critical values of Landau-Ginzburg potentials, we observe slightly 
different properties from the Landau-Ginzburg theory of Fano 
varieties such as the projective spaces \cite{Du2,Gu}. 

\vskip0.5cm

Let $X$ be a quintic Calabi-Yau hypersurface in 
the projective space of dimension four, and $Y$ be its mirror 
hypersurface\cite{CdOGP}. In the application to mirror symmetry,  
the mirror hypersurface family appears as a one-parameter family $\{Y_x\}$ 
of Calabi-Yau hypersurfaces in a toric variety. 
In the torus $(\mathbf C^*)^4$, which 
is dense in the toric variety, we may describe $Y_x$ explicitly by  
first considering $Y_a$ as a hypersurface $f(y,a)=0$ with 
\begin{equation}
f(y,a)=y_1+y_2+y_3+y_4+\frac{1}{y_1y_2y_3y_4} + a \; .
\mylabel{eqn:LG-f}
\end{equation}
Since the two hypersurfaces $f(y,a)=0$ and 
$f( y, \alpha a)=0$ $(\alpha^5=1)$ are isomorphic under 
$y \rightarrow \alpha y$, we identify $Y_a$ with $Y_{\alpha a}$ 
introducing the parameter $x=-\frac{1}{a^5}$ to define the 
mirror family $\{ Y_x \}$. 

Now let us take a cycle $\gamma \in H_3(Y_{x_0},\mathbf Z)$ and 
define the following period integral 
\begin{equation}
\Pi_\gamma(x)=
\frac{1}{(2\pi i)^4}
\int_\gamma Res_{f(y,a)=0}\left( \frac{a}{f(y,a)} 
\prod_{k=1}^4\frac{dy_k}{y_k} \right) \;, 
\mylabel{eqn:Pi[x]}
\end{equation}
for the family. It is easy to derive the following 
Picard-Fuchs differential equation,
\begin{equation}
\{ \theta_x^4 -5 x 
(5\theta_x+4)
(5\theta_x+3)
(5\theta_x+2)
(5\theta_x+1) \} \Pi_\gamma(x) =0\;.
\mylabel{eqn:PF-Y5}
\end{equation}
Our interest in this note is the following Laplace transform of 
this Picard-Fuchs differential equation, 
\begin{equation}
\{ z(\theta_z+1)^3+5(5\theta_z+1)(5\theta_z+2)(5\theta_z+3)(5\theta_z+4) \}
\hat \Pi_\Gamma(z) =0 \;, 
\mylabel{eqn:Qdiff-Y5}
\end{equation}
where $\hat\Pi_\Gamma(z)$ is the (formal) Laplace transform of the 
the period integral. When the cycle $\gamma$ is a vanishing cycle,  
the Laplace transform $\hat\Pi_\Gamma(z)$ may be written more precicely by 
\begin{equation}
\begin{aligned}
\int_{x_c}^\infty e^{-zx} \Pi_\gamma(x) dx 
&
= \frac{1}{(2\pi i)^4} \int_{x_c}^\infty e^{-zx} 
\int_\gamma Res\left( 
\frac{1}{W(y)-x} \prod_k dy_k \right) dx  \\
& 
= \frac{1}{(2\pi i)^4} 
\int_\Gamma e^{-z W(y)} dy_1 dy_2 dy_3 dy_4 \;\;,
\end{aligned}
\mylabel{eqn:Osci-Int}
\end{equation}
where we set 
\begin{equation}
W(y)=y_1y_2y_3y_4(y_1+y_2+y_3+y_4+1) \;\;,
\mylabel{eqn:LGpot-W}
\end{equation} 
and $\Gamma$ represents the Lefschetz thimble made by the vanishing 
cycle $\gamma$ which vanishes at the critical value $x_c$ of $W(y)$ 
\cite{AGV}. 
$W(y)$ is called the Landau-Ginzburg potential, and $\hat \Pi_\Gamma(z)$ 
above is the oscillatory integral of the potential.  
It is easy to see that the Laplace transform (\ref{eqn:Qdiff-Y5})  
follows from (\ref{eqn:PF-Y5}) by simple replacements 
$\theta_x \rightarrow -\theta_z-1$ and 
$\frac{\partial \; }{\partial x} \rightarrow z$ 
after a division by $x$. 

The asymptotics of the oscillatory integral $\hat \Pi_\Gamma(z)$ 
when $z\rightarrow \infty$ are determined by the critical values of 
the potential function $W(y)$. We find that there are two critical 
values $x_c=0$ and $x_c=\frac{1}{5^5}$, and observe a (continuous) 
degeneration of the critical points in $W^{-1}(0)$ yet only an 
isolated critical point in $W^{-1}(\frac{1}{5^5})$. 
The degeneration at the critical value $0$ is characterized by 
maximal unipotent monodromy there \cite{Mo}. 
In our case, this degeneration results in a certain distinguished 
property of the solutions about the irregularity at $z=\infty$ of 
$(\ref{eqn:Qdiff-Y5})$. The purpose and the main 
results of this note involve analysis of the asymptotics when 
$z\rightarrow\infty$ under such a
degeneration for the mirror quintic and the similar hypersurfaces 
studied in \cite{KT}. 

Also, more generally, we introduce two different kinds of 
Laplace transforms of the Gel'fand-Kapranov-Zelevinski (GKZ) system, 
$\widehat{\rm GKZ}^*$ and $\widehat{\rm GKZ}_\nu$, and derive the 
Laplace transform (\ref{eqn:Qdiff-Y5}) 
from these general systems of differential equations.

\vskip0.3cm

This note is organized as follows: In Section 2, we briefly summarize 
the standard definition of the Stokes lines, and then present our analysis 
of the differential equation (\ref{eqn:Qdiff-Y5}). In Section 3, we will 
introduce the Laplace transforms $\widehat{\rm GKZ}^*$ and 
$\widehat{\rm GKZ}_\nu$ of the relevant GKZ system, and derive 
(\ref{eqn:Qdiff-Y5}) from a $\widehat{\rm GKZ}_\nu$ system. Conclusions and 
discussions are presented in Section 4, where we interpret the Stokes 
matrices from mirror symmetry. In Appendix A, we summarize similar 
results for other Calabi-Yau hypersurfaces in weighted 
projective spaces. 

\vskip0.5cm

\noindent
{\bf Acknowledgments:}
C.F.D. would like to thank the Graduate School of Mathematical Sciences, 
University of Tokyo and the Perimeter Institute, Waterloo for 
supporting his visits
during the critical early stages of this project.
S.H. would like to thank for the organizers of the research programs 
``Mathematical Structures in String Theory 2005'' at KITP, Santa Barbara, 
and ``Modular Form and String Duality'' at BIRS, Banff (Jun.~2006) for 
providing him nice research environments where this work has 
progressed a lot. 
We thank A. Corti, J. Morgan, and M.-H. Saito for valuable discussions.
C.F.D. is supported in part by a Royalty Research Fund Scholar Award from
the Office of Research, University of Washington. 
S.H. is supported in part by the Grant-in Aid for Scientific Research 
C18540014.

\vskip1cm
\section{{\bf Degeneration and Stokes matrices}}

{\bf (2-1) Stokes matrices -- non-degenerate case. } 
Here we summarize the description of the asymptotic 
solutions about an irregular singular point and the definition of 
Stokes matrices in general. This is to set our notations and also 
to contrast the general cases with our degenerate cases.  

\vskip0.5cm
For simplicity, let us consider a linear differential equation of order $n$  
which has a regular singularity at $z=0$ and an irregular singularity 
at $z=\infty$. For example, the differential equation which describes 
the (small) quantum cohomology of the projective space $\mathbf P^{n-1}$ is 
of this type (see (\ref{eqn:Qdiff-IP4}) below). 
About $z=0$ we can construct power series 
solutions of the differential equation with infinite radius of convergence;
we denote them by $\phi_1(z),\cdots,\phi_n(z)$. For the 
solutions about $z=\infty$, we make the following ansatz for 
the (formal) solution,
\begin{equation}
g(z)=e^{\Lambda(z)} z^{-r}p(z)
    =e^{\Lambda(z)} z^{-r}\left( 1+\frac{C_1}{z}+\frac{C_2}{z^2}+\cdots \right), 
\mylabel{eqn:asympt-sol}
\end{equation}
where $\Lambda(z)$ is a polynomial. According to the general theory, 
by substituting this ansatz 
into the differential equation, we can determine the 
polynomial $\Lambda(z)$ and also the constants $r, C_1,C_2, \cdots $. 
In this paper we will call the differential equation {\it non-degenerate} 
if we can find $n$ formal solutions of the above form with linearly 
independent polynomials $\Lambda_k(z)$ $(k=1,\cdots,n)$, and 
otherwise {\it degenerate}. 

\vskip0.5cm
It is a standard result that even though the series in (\ref{eqn:asympt-sol}) 
is a formal one, it still provides a way to describe the asymptotics 
of the holomorphic solutions $\phi_k(z)$ near $z=\infty$. 
Suppose the differential equation is non-degenerate, and write $n$ 
independent formal solutions $g_k(z)=e^{\Lambda_k(x)}z^{-c_k}p_k(z) 
\; (k=1,\cdots,n)$. When $z \rightarrow \infty$, the relative magnitude  
of $|g_k(z)|$ and $|g_l(z)|$ is governed by the sign of 
${\rm Re}(\Lambda_k(z) - \Lambda_l(z))$. Namely, $|g_k(z)/g_l(z)|>1$ 
(respectively, $<1$) for the asymptotic region defined by 
${\rm Re}(\Lambda_k(z) - \Lambda_l(z))>0$ $(<0)$. Note that the asymptotic 
region consists of $d$ connected components, where $d$ is  
the degree of the polynomial $\Lambda_k(z) - \Lambda_l(z)$. Each connected 
component in $|z|\gg1$ of ${\rm Re}(\Lambda_k(z) - \Lambda_l(z))>0$ (respectively $<0$) 
is called a {\it positive (negative) angular region} of 
${\rm Re}(\Lambda_k(z) - \Lambda_l(z))$.  There are $2d$ angular 
regions in total. We will find it convenient, however, 
to distinguish the angular region 
written by $\theta < arg(z) < \theta'$ $(|z|\gg1)$ 
with its $2\pi$ rotation, i.e. $\theta + 2\pi < arg(z) < \theta' + 2 \pi$ 
$(|z|\gg1)$. 
The {\it Stokes lines} are defined to be the half lines from 
the origin which separate these asymptotic angular regions. 

\vskip0.5cm
Since, for the differential equations considered in this note, the polynomial 
$\Lambda_k(z)$ is always linear, we assume hereafter that 
$\Lambda_k(z)=c_k z$ with some constant $c_k$. In this case, the line 
${\rm Re}(\Lambda_k(z)-\Lambda_l(z))=0$ simply 
divides the complex $z$-plane into two connected components, the positive 
and negative angular regions. We split the line 
${\rm Re}(\Lambda_k(z)-\Lambda_l(z))=0$ 
into the two Stokes rays, $R_{kl}$ and $R_{lk}$, so that  
the counter-clockwise rotation of the half line $R_{ij}$ belongs to 
the positive angular region of ${\rm Re}(\Lambda_i(z)-\Lambda_j(z))$  
(see Fig.~1). The notation $\overline R_{kl}$ for $R_{lk}$ will be used 
in the following. 


\vskip0.5cm
\centerline{\epsfxsize 9.5truecm\epsfbox{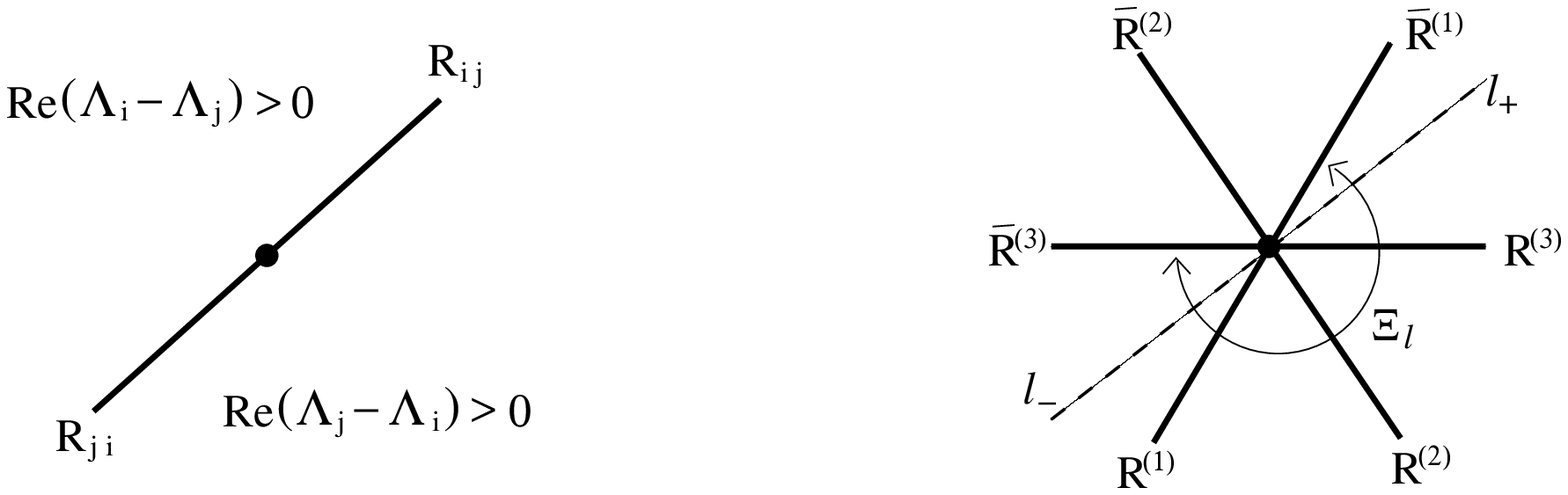}}

{\leftskip1cm \rightskip1cm 
\noindent 
{{\bf Fig.1.}  The convention of the Stokes ray $R_{ij}$ explained in the 
text (left). An example of Stokes lines and the oriented line 
$l=(l_-,l_+)$(right). $\Xi_l$ shows the angular region determined by 
the line $l$. \par} }


\vskip0.5cm
We define
$$
\Xi (\underline\theta,\bar\theta):=\{ \;z  \; | \; 
\underline\theta < arg(z) < \bar\theta \; \}, 
$$
and say that $\Xi(\underline\theta,\overline\theta)$ is a 
{\it proper angular region} of $\Lambda_k(z)$ if for each $j\not=k$ 
there is at most one positive angular region of 
${\rm Re}(\Lambda_j(z)-\Lambda_k(z))$ which has non-empty intersection with 
$\Xi(\underline\theta,\overline\theta)$. The following theorem is standard 
(see e.g. [CL]):

\vskip0.5cm
\begin{theorem} \mylabel{thm:thmRegSol} Fix the formal solutions $g_k(z) \; 
(k=1,\cdots,n)$ as above. If $\Xi(\underline\theta,\overline\theta)$ is 
a proper angular region of 
for all $\Lambda_k(z)$ $(k=1,\cdots,n)$, then there are $n$ independent  
holomorphic solutions $f_i(z)$ $(i=1,\cdots,n)$ 
whose asymptotic expansions in $\Xi(\underline\theta,\overline\theta)$  
are given by 
$$
f_i(z) \sim g_i(z)=e^{\Lambda_i(z)}z^{-r_i}\left(1
+\frac{C_{i1}}{z}+\frac{C_{i2}}{z^2} + \cdots \right)   \;\;\; 
(z\rightarrow \infty, i=1,\cdots,n). 
$$
Moreover such solutions $f_i(z)$ are unique if $\Xi(\underline\theta, 
\overline\theta)$ is not contained in  any negative angular region 
$Re(\Lambda_i(z)-\Lambda_j(z))$ $(1\leq i,j \leq n)$. 
\end{theorem}

\vskip0.5cm
To make use of this theorem, it is convenient to set notation for 
an angular region on which the solutions $f_i(z)$ above  
are determined uniquely. Let us consider a line $l$ passing through 
the origin and located in a general position to the Stokes lines.  
We introduce the orientation of $l$ by writing it as a union of 
the half lines from the origin $l_-,l_+$, i.e. $l=(l_-,l_+)$. 
We define the angular region $\Xi_l$ to be the oriented line given by
\begin{equation*}
\Xi_l=\Xi(\theta_{l_-}-\epsilon_-, \theta_{l_+}+\epsilon_+) ,
\end{equation*}
where $\theta_{l_-}$ represents $arg(z)$ for $z \in l_-$ and 
$\theta_{l_+}=\theta_{l_-}+\pi$. The angle $\epsilon_- >0$ 
is determined so that $\theta_{l_-}-\epsilon_-$ represents the angle  
of the nearest Stokes ray with smaller angle than $l_-$. Similarly 
for $\epsilon_+ >0$ we use the nearest Stokes ray with larger angle 
than $l_+$. More precisely, $\Xi_l$ is defined only up to $2\pi$ rotations, 
however this ambiguity will be fixed in each case when 
it will become necessary to do so. 
Also the naive notaions $\Xi_{l+\pi}$, $\Xi_{l-\pi}$, ... are useful 
if we understand $l+\pi$ (respectively, $l-\pi$) to be the counter-clockwise 
(clockwise) rotation by $\pi$ of the oriented line $l=(l_-,l_+)$. 
It is an easy exercise to see that the angular region $\Xi_l$ is 
proper for all $\Lambda_k(z)$ and, in addition, 
is not contained in any negative angular region of 
${\rm Re}(\Lambda_i(z)-\Lambda_j(z))$ $(1\leq i,j \leq n)$.

\vskip0.5cm
In this note, we adopt the convention of arranging the formal 
solutions $g_k(z)$ and the holomorphic solutions $f_k(z)$ into row vectors:
\begin{equation*}
\Psi(z)=(g_1(z),g_2(z),\cdots,g_n(z)) \;\;,\;\;
\Phi(z)=(f_1(z),f_2(z),\cdots,f_n(z)) .
\end{equation*}
The following theorem is also standard (see e.g. \cite{CL}):

\vskip0.5cm
\begin{theorem} 
Fix the formal solutions $g_k(z)$, and fix the angular 
region $\Xi_l$ with an oriented line $l=(l_-,l_+)$. Let $l'=(l_-',l_+')$ 
be the oriented line obtained by rotating $l_-$ in a counter-clockwise 
direction and 
passing $l_-$ through the nearest Stokes ray $R$. Denote the corresponding 
unique solutions by $\Phi_l(z)$ and $\Phi_{l'}(z)$, respectively. 
If the Stokes ray is written $R=R_{ij}$, then these two solutions are 
related by 
$$
\Phi_{l'}(z)=\Phi_l(z) K_{R_{ij}}  \;\;, \;\;  z \in \Xi_l \cap \Xi_{l'} \;,
$$
with 
\begin{equation*}
K_{R_{ij}}=I_n + c_{R_{ij}} E_{ji} \;,
\end{equation*}
where $c_{R_{ij}}$ is a constant, $I_n$ is the identity matrix, and 
$E_{ji}$ is the matrix with non-vanishing entry $1$ at the $ji$ 
position and otherwise $0$. Similarly, let $l''$ be the oriented
line obtained by clockwise movement of 
$l_+$ passing through the nearest Stokes ray $R'=R_{jk}$, and 
$\Phi_{l''}(z)$ be the corresponding holomorphic solution on $\Xi_{l''}$, 
then 
\begin{equation*}
\Phi_{l''}(z) = \Phi_l(z) K_{R_{jk}}^T \;\;,\;\; z \in \Xi_l \cap \Xi_{l''},  
\end{equation*}
where the superscript $T$ represents the transpose.
\end{theorem}

\vskip0.5cm
Note that the constant $c_{R}$ above depends only on the Stokes ray $R$, 
and also the following relation is implicit in the above theorem;
\begin{equation}
K_{R_{jk}}^T = K_{R_{kj}}^{-1} \;\; (= K_{\overline R_{jk} }^{-1})\;,
\mylabel{eqn:KT}
\end{equation}
for all Stokes rays. Successive applications of the above theorem 
give
\begin{equation*}
\begin{aligned}
\Phi_{l+\pi}(z)&=\Phi_{l}(z) \, K_{R^{(1)}}K_{R^{(2)}} \cdots K_{R^{(n)}} 
\;\;\;\;\; (z \in \Xi_+), \cr
\Phi_{l-\pi}(z)&=\Phi_{l}(z) \, K_{R^{(n)}}^T 
                              K_{R^{(n-1)}}^T \cdots 
                              K_{R^{(1)}}^T 
\;\; (z \in \Xi_-),
\end{aligned}
\end{equation*}
where $R^{(1)},\cdots,R^{(n)}$ are the Stokes rays in the order through which  
$l_-$ passes during the rotation by $\pi$, and $\Xi_{\pm}$ 
are the intersections of the successive angular regions in the process.
In the literature, 
$\Phi_{l}(z), \Phi_{l+\pi}(z)$ and $\Phi_{l-\pi}(z)$ are often denoted 
by $\Phi^{right}(z), \Phi^{left}(z)$ and $\Phi^{left'}(z)$, respectively, 
and the above linear relations are written 
$\Phi^{left}(z)=\Phi^{right}(z) S$ and 
$\Phi^{left'}(z)=\Phi^{right}(z) S_-$ with 
\begin{equation*}
S=K_{R^{(1)}}K_{R^{(2)}} \cdots K_{R^{(n)}} \;\;,\;\;  
S_-=K_{R^{(n)}}^T K_{R^{(n-1)}}^T \cdots K_{R^{(1)}}^T \;\;. 
\end{equation*}
The matrices $S, S_-$ are the so-called {\it Stokes matrices}. 
By definition, the Stokes matrix $S_-$ satisfies $S_-=S^T$, and 
the following relation holds:
\begin{equation*}
\Phi^{left}(z)=\Phi^{left'}(e^{-2\pi i}z) M_0 \;\;,\;\;
M_0=(S^{-1})^T S \;\;, 
\mylabel{eqn:M0}
\end{equation*}
where $M_0$ represents the monodromy about the regular singularity at $z=0$.

\vskip1cm
{\bf (2-2) Stokes matrices -- degenerate case. } 
Now let us turn our attention to our differential equation 
(\ref{eqn:Qdiff-Y5}) whose irregular singularity is degenerate according 
to the definition given in (2-1). 

\vskip0.5cm
Let us first determine the formal solutions $g_k(z)$ near the 
irregular singularity at $z=\infty$. To describe the solutions, we first 
define the following hypergeometric series:
$$
G(z,\rho)=\sum_{n=0}^\infty \frac{\Gamma(1+5(n+\rho))}{\Gamma(n+\rho+1)^4}
\left( \frac{1}{z} \right)^{n+\rho+1} \;\;.
$$
Then it is straightforward to obtain (and define) the following formal 
solutions.

\vskip0.5cm
\begin{proposition} 
The differential equation (\ref{eqn:Qdiff-Y5}) is degenerate, and the 
following $g_k(z)$ are the formal solutions at the irregular singularity 
$z=\infty$:
\begin{equation}
\begin{aligned}
& 
g_1(z)=G(z,0) \;,\;\; 
g_2(z)=\pd_{\tilde\rho}G(z,\rho)\Big|_{\rho=0} ,\;\;
g_3(z)=\big\{ \frac{5}{2} \pd_{\tilde\rho}^2+{\tt a} \pd_{\tilde\rho} \big\} 
G(z,\rho)\Big|_{\rho=0} \\
&
g_4(z)=\frac{5^5\sqrt{5}}{2\pi i} e^{-\frac{z}{5^5}} \frac{1}{z^2} 
\left(1-\frac{4375}{z}+\frac{32031250}{z^2}-\cdots \right)\;, \\
\end{aligned}
\mylabel{eqn:formal-g-quintic}
\end{equation}
where $\pd_{\tilde\rho}=\frac{1}{2\pi i} \frac{\pd \;}{\pd \rho}$ and 
${\tt a}$ is a parameter to be discussed later. 
\end{proposition}

\vskip0.5cm
\centerline{\epsfxsize 3truecm\epsfbox{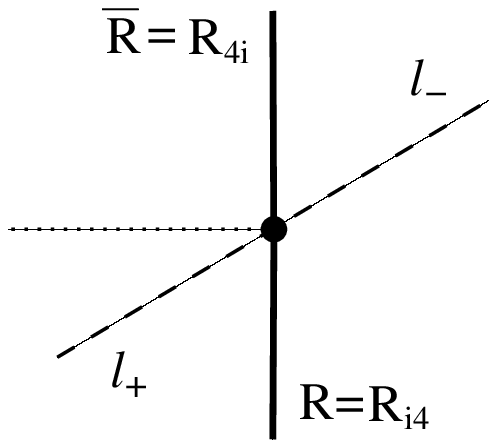}}
{\leftskip1cm \rightskip1cm 
\noindent 
{{\bf Fig.2.} The Stokes rays of (\ref{eqn:formal-g-quintic}) 
and the cut of the logarithms (the dotted line) 
for Calabi-Yau hypersurfaces $Y(d)$ $(d=5,6,8,10)$.  
 \par} }
\vskip0.5cm

As we see above $\Lambda_1(z)=\Lambda_2(z)=\Lambda_3(z)=0$ and $\Lambda_4(z)
=-\frac{z}{5^5}$, hence the differential equation is degenerate. 
The Stokes rays $\overline R=R_{4i}$ and $R=R_{i4}$ ($i=1,2,3$), and 
also the oriented line $l=(l_-,l_+)$ are defined as before, although their 
configulation is slightly different (see Fig.~2). 
For the oriented line $l$, we fix the angular region $\Xi_l$ by 
\begin{equation}
\Xi_l=\Xi(-\frac{\pi}{2},\frac{3\pi}{2})
=\Xi(-\frac{3\pi}{2},\frac{3\pi}{2})\cap \Xi(-\frac{\pi}{2},\frac{5\pi}{2}) \;,
\mylabel{eqn:angR-quintic}
\end{equation}
and write $\Xi_{l+\pi m}$ for the $\pi m$ rotations of $\Xi_l$, i.e. 
$\Xi_{l+m\pi}=\Xi(-\frac{\pi}{2}+m\pi,\frac{3\pi}{2}+m\pi)$.
Note that the angular region $\Xi(-\frac{3\pi}{2},\frac{3\pi}{2})$ 
is proper for $\Lambda_4(z)$, and $\Xi(-\frac{\pi}{2},\frac{5\pi}{2})$ 
is proper for $\Lambda_i(z) \; (i=1,2,3)$.

\vskip0.5cm

The solutions about the regular singularity $z=0$ are given by
\begin{equation*}
\phi_{k}(z)=\sum_{n=0}^\infty 
\frac{\Gamma(1+n+\rho_k)^3\Gamma(2+n+\rho_k)}{\Gamma(6+5(n+\rho_k))}
(-1)^n z^{n+\rho_k} \;\;,
\quad 
\end{equation*}
where $\rho_k=-\frac{k}{5}$ ($k=1,\cdots,4$). 
We may write these solutions in terms of Barnes hypergeometric series
\begin{equation}
\begin{aligned}
&
\phi_1(z)=
\frac{d_1}{z^{\frac{1}{5}}} 
\,_3F_3(\frac{4}{5}^3;\frac{6}{5},\frac{7}{5},\frac{8}{5},
-\frac{z}{5^5}) ,\;\;
\phi_2(z)=
\frac{d_2}{z^{\frac{2}{5}}}
\,_3F_3(\frac{3}{5}^3;\frac{4}{5},\frac{6}{5},\frac{7}{5},
-\frac{z}{5^5}) \;\;,\;\; \\
&
\phi_3(z)=
\frac{d_3}{z^{\frac{3}{5}}} 
\,_3F_3(\frac{2}{5}^3;\frac{3}{5},\frac{4}{5},\frac{6}{5},
-\frac{z}{5^5}) ,\;\;
\phi_4(z)=
\frac{d_4}{z^{\frac{4}{5}}}
\,_3F_3(\frac{1}{5}^3;\frac{2}{5},\frac{3}{5},\frac{4}{5},
-\frac{z}{5^5}) \;\;,\;\; \\
\end{aligned}
\end{equation}
where the constants $d_1,d_2,d_3,d_4$ are given, respectively, by 
$\frac{\Gamma(\frac{4}{5})^4}{30}$, 
$\frac{\Gamma(\frac{3}{5})^4}{10}$, 
$\frac{\Gamma(\frac{2}{5})^4}{5}$, and
$\frac{\Gamma(\frac{1}{5})^4}{5}$.

To make our definitons parallel to those of the hypergeometric 
series which appear in \cite{CdOGP}, we set
\begin{equation}
\hat w^{\infty}_m(z)=  (1-\alpha)
\sum_{k=1}^4 (1-\alpha^k)^3\alpha^{km+\frac{1}{2}(k-1)}\,\phi_{5-k}(z) \;\;\;
(m=0,1,\cdots,4),
\end{equation}
with $\alpha=e^{2\pi i/5}$. Note that, by definition, we have 
$$
\hat w^\infty_0(z)+
\hat w^\infty_1(z)+ \cdots + \hat w^\infty_4(z) = 0 \;\;.
$$
Now we define the holomorphic solutions $f_k(z)$ about $z=0$ by 
\begin{equation}
(\;f_1(z)\,f_2(z)\,f_3(z)\,f_4(z)\;) = c_N \, 
(\;\hat w^\infty_0(z)\, \hat w^\infty_1(z)\, 
   \hat w^\infty_2(z)\, \hat w^\infty_4(z)\;) \,N \;\;,
\mylabel{eqn:fk-phi-quintic}
\end{equation}
where $c_N$ is a normalization constant and $N$ is the matrix given 
in Table 1 below.  The asymptotics of the solutions $f_k(z)$ are 
given by the following result:

\vskip0.5cm
\begin{proposition} \mylabel{thm:fk-gk}
There is a constant $c_N$ for which the following asymptotics hold 
as $z \rightarrow \infty$:
\begin{equation}
\begin{aligned}
& f_k(z) \;\sim\; g_k(z) \;\; 
(z \in \Xi(-\frac{\pi}{2},\frac{5\pi}{2}), \; k=1,2,3)  \;,\\
& f_4(z) \;\sim\; g_4(z) \;\; (z \in \Xi(-\frac{3\pi}{2},\frac{3\pi}{2}) ) .\\
\end{aligned}
\mylabel{eqn:fg-asym-quintic}
\end{equation}
Hence for the angular 
region $\Xi_l=
\Xi(-\frac{\pi}{2},\frac{5\pi}{2}) \cap \Xi(-\frac{3\pi}{2},\frac{3\pi}{2})$ 
we have 
$$
\Phi_l(z)=(\, f_1(z)\, f_2(z) \, f_3(z) \, f_4(z) \,) \; \sim \; 
(\, g_1(z) \, g_2(z) \, g_3(z) \, g_4(z) \,) 
\;\; (z \rightarrow \infty).
$$
\end{proposition}

\vskip0.5cm 
Our results above rely on numerical calculations of 
order $10^{-50}$ by {\it Maple}. In the numerical calculation, the 
constant $c_N$ is $-0.00342934921...$ \, . 

\vskip0.5cm

To evaluate the Stokes matrices, we arrange the asymptotic 
solutions $g_k(z)$ as 
\begin{equation*}
\Psi(z)=(g_1(z)\,g_2(z) \, g_3(z) \, g_4(z)) \;\;,
\end{equation*}
and note that $\Psi(z)$ has the following monodromy due to the 
logarithms in $g_k(z)$:
\begin{equation*}
\Psi(e^{-2\pi i}z)=\Psi(z) \, L \;\;,\;\; 
L=\left( \begin{matrix} 
1 & 1 & \frac{5}{2}+{\tt a} & 0 \\
0 & 1 & 5 & 0 \\
0 & 0 & 1 & 0 \\
0 & 0 & 0 & 1 \\ \end{matrix} \right) \;\;.
\end{equation*}
Now define the holomorphic solutions $\Phi_{l+m\pi}(z)$ about $z=0$ by 
\begin{equation*}
\Phi_{l+m\pi}(z)=
\begin{cases}
(\, f_1(z_{m-1})\,f_2(z_{m-1})\,f_3(z_{m-1})\,
   f_4(e^{-2\pi i}z_{m-1})\,) \,
L^{-\frac{m-1}{2}} 
& m: odd  \\
(\,f_1(z_m)\,f_2(z_m)\,f_3(z_m)\,f_4(z_m)\,)\,L^{-\frac{m}{2}}
& m: even \\ 
\end{cases} 
\end{equation*}
where we set the notation $z_s=e^{-s\pi i}z$ for integer $s$.

\vskip0.5cm
\begin{proposition} \mylabel{thm:mainI} 
The holomorphic solutions $\Phi_{l+m\pi}(z) \; (m \in \mathbf Z)$ have 
the asymptotics: 
\begin{equation}
\Phi_{l+m\pi}(z) \;\sim\; \Psi(z) \quad (z \in \Xi_{l+m\pi}, \; z\rightarrow \infty).
\mylabel{eqn:asympt-lm}
\end{equation}
They satisty the following linear relations 
\begin{equation}
\begin{aligned}
\Phi_{l+(2k+1)\pi}(z)&=\Phi_{l+2k \pi}(z) \, L^{k}K_{\overline R} L^{-k}  \;,\\
\Phi_{l+(2k+2)\pi}(z)&=\Phi_{l+(2k+1) \pi}(z) \, L^{k+1} K_{R} L^{-k-1} \quad, 
\end{aligned}
\mylabel{eqn:connectKL-quintic}
\end{equation}
on the respective angular regions $\Xi_{l+(2k+1)\pi} \cap \Xi_{l+2k \pi}$ 
and $\Xi_{l+(2k+2)\pi} \cap \Xi_{l+(2k+1) \pi}$, with the matrices 
$K_{\overline R}, K_R$ defined by 
\begin{equation*}
K_{\overline R}=\left( 
               \begin{matrix} 1 & 0  & 0 &-5\\
                             0 & 1 & 0 &-\frac{5}{2}+{\tt a} \\
                             0 & 0 & 1 &-1\\ 
                             0 & 0 & 0 & 1\end{matrix}\right) \;\;,\;\;
K_{R}=\left( 
                \begin{matrix} 1 & 0 & 0 & 0 \\
                             0 & 1 & 0 & 0 \\
                             0 & 0 & 1 & 0 \\
                             1 & 0 & 0 & 1 \\ \end{matrix}\right) \;\;.\;\;
\end{equation*}
\end{proposition}

\vskip0.3cm
\begin{proof}
The first half of the claim follows from Proposition 
\ref{thm:fk-gk} and the definition of $\Phi_{l+m\pi}(z)$. For example, 
when $m=2k+1$, we have for 
$z \in \Xi_{l+(2k+1)\pi}
=\Xi(-\frac{\pi}{2}+(2k+1)\pi,\frac{3\pi}{2}+(2k+1)\pi)$,
\begin{equation*}
\begin{aligned}
\Phi_{l+(2k+1)\pi}(z)
&=(\, f_1(e^{-2k\pi i}z) \; f_2(e^{-2k\pi i}z) \; 
      f_3(e^{-2k\pi i}z) \; f_4(e^{-(2k+2)\pi i}z)\,)
\,L^{-k} \\
&\sim (\, g_1(e^{-2k\pi i}z) \; g_2(e^{-2k\pi i}z) 
      \; g_3(e^{-2k\pi i}z) \; g_4(e^{-(2k+2)\pi i}z\,)
\,L^{-k}  \\
&=(\, g_1(z) \; g_2(z) \; g_3(z) \; g_4(z) \,)\;, \\
\end{aligned}
\end{equation*}
where we have used the relation (\ref{eqn:fg-asym-quintic}) 
in the second line, and also $g_4(e^{2\pi i}z)=g_4(z)$ in the third line. 
Likewise we can verify the claimed asymptotics for $m=2k$.

For the second half of the claim, let us note 
\begin{equation*}
\phi_k(e^{-2\pi i}z)=\alpha \phi_k(z) \;\;\;(k=1,2,3,4). 
\end{equation*}
Then, after some linear algebra using the definition 
(\ref{eqn:fk-phi-quintic}), it is straightforward to verify 
the following identities, 
\begin{equation*}
\begin{aligned}
& (\;f_1(z)\;f_2(z)\;f_3(z)\;f_4(z')\;)=
(\;f_1(z) \;f_2(z) \; f_3(z) \; f_4(z) \;) \, K_{\overline R} \;,\\
& (\;f_1(z')\;f_2(z')\;
     f_3(z')\;f_4(z') \;)=
  (\;f_1(z)\;f_2(z)\;f_3(z)\;f_4(z')\;) \, L K_{R} \;\;,\\
\end{aligned}
\end{equation*}
with $z'=e^{-2\pi i}z$ and the given matrices $K_{\overline R}$ and $K_R$. 
These relations are nothing but $\Phi_{l+\pi}(z)=\Phi_l(z) K_{\overline R}$ 
and $\Phi_{l+2\pi}(z)L=\Phi_{l+\pi}(z) LK_R$. The general relations 
(\ref{eqn:connectKL-quintic}) follow from the above identities applied to 
$z'=e^{-2k\pi i}z$, i.e.,   
\begin{equation*}
\begin{aligned}
\Phi_{l+(2k+1)\pi}(z)&=
(\,f_1(z')\;f_2(z')\;f_3(z')\;f_4(e^{-2\pi i}z')\,)\,L^{-k} \\
&= (\,f_1(z')\;f_2(z')\;f_3(z')\;f_4(z')\,)\,K_{\overline R} L^{-k} 
= \Phi_{l+2k \pi}(z) L^{k}K_{\overline R} L^{-k} , \\
\end{aligned}
\end{equation*}
and a similar formula for the second relation. 
\end{proof}

\vskip0.5cm
\begin{remark} 
(1) At first sight, the relation (\ref{eqn:connectKL-quintic}) 
looks incompatible 
with the asymptotic behavior (\ref{eqn:asympt-lm}) when we take 
$z \rightarrow \infty$. However this is not the case due to 
the triangular nature of the matrices $K_{\overline R}, K_R, L$ 
and the asymptotic properties of $g_k(z)$ .

(2) By the relations (\ref{eqn:connectKL-quintic}), we have 
$\Phi_{l+\pi}(z)=\Phi_l(z) K_{\overline R}$ for $z \in \Xi_l\cap\Xi_{l+\pi}$ 
and $\Phi_{l-\pi}(z) = \Phi_l(z) K_R^{-1}$ 
for $z \in \Xi_l\cap\Xi_{l-\pi}$. From these we read the Stokes matrices as  
$$
S=K_{\overline R} \;\;,\;\;  
S_-=K_{R}^{-1}.
$$
Here one should note that the relation (\ref{eqn:KT}), which guarantees 
the symmetry of exchanging $l_-$ and $l_+$ in the argument there, is no longer 
valid for our degenerate case. 
Also the Stokes matrices $S$ and $S_-$ are not constant, 
but must be conjugated by $L$ in a way that depends on the sheets where 
the holomorphic solutions  $\Phi_{l+m \pi}(z)$ are defined. Of course 
the latter property simply arises from the logarithmic behavior of the 
formal series $\Psi(z) =(g_1(z) \,g_2(z)\, g_3(z)\,g_4(z))$. 
\end{remark}

\vskip1cm

{\bf (2-3) Stokes matrices and monodromy matrices.} 
There are simple relations between the Stokes 
matrices $K_R$, $K_{\overline R}$ and the monodromy matrices of the 
period integrals $\Pi_{\gamma}(x)$. Let us recall that the Picard-Fuchs 
differential equation (\ref{eqn:PF-Y5}) has regular singularities at 
$x=0, \frac{1}{5^5}, \infty$, and that the integral symplectic bases of 
its solutions, first obtained in \cite{CdOGP}, can be written 
in a concise way as follows (see \cite{Ho1,Ho2} for some 
generalizations): 

\vskip0.5cm
\begin{proposition} \mylabel{thm:sympPi} 
About the regular singularity at $x=0$ ($|x|<\frac{1}{5^5}$), 
the integral and symplectic bases of the solutions of the 
Picard-Fuchs equation (\ref{eqn:PF-Y5}) are given by 
\begin{equation}
\begin{aligned}
&\Pi_{\gamma_1}(x) = w(x,0) \;\;,\;\;  
&\Pi_{\gamma_3}(x)=\big\{ \frac{5}{2}\pd_{\tilde\rho}^2 + 
               {\tt a} \pd_{\tilde\rho} \big\} w(x,\rho)\big|_{\rho=0} \;\;, 
   \hskip0.5cm \\
&\Pi_{\gamma_2}(x)=\pd_{\tilde\rho} w(x,\rho)\big|_{\rho=0} \;\;,\;\;
&\Pi_{\gamma_4}(x)=\big\{ - \frac{5}{3!}\pd_{\tilde\rho}^3 
      - \frac{50}{12} \pd_{\tilde\rho} \big\} w(x,\rho)\big|_{\rho=0} \;\;, 
   \\
\end{aligned}
\mylabel{eqn:periods-quintic}
\end{equation}
where $w(x,\rho)=\sum_{n=0}^\infty 
\frac{\Gamma(1+5(n+\rho))}{\Gamma(1+n+\rho)^5} x^{n+\rho}$. The solutions
$\Pi_{\gamma_k}(x)$ and $\Pi_{\gamma_{5-k}}(x)$ $(k=1,2)$ are 
symplectic dual to each other and the constant 
${\tt a} \in \mathbf Z+\frac{1}{2}$ is an arbitrary half (odd) integer. 
\end{proposition}

\vskip0.5cm 
It is known that the cycle $\gamma_1 \approx T^3$ is a torus cycle 
coming from the ambient toric variety \cite{CdOGP}\cite{Ba2}. The 
dual cycle $\gamma_4=S^3$ is a vanishing cycle that 
appears at $x=\frac{1}{5^5}$ \cite{CdOGP}. In our form of the 
Laplace transform (\ref{eqn:Osci-Int}), we can identify this cycle 
as the vanishing cycle associated to the critical value 
$x_c=\frac{1}{5^5}$ of the LG potential $W(y)$ in (\ref{eqn:LGpot-W}). 
In fact, after an analytic continuation of $\Pi_{\gamma_4}(x)$, we  
obtain a convergent powerseries in $t=x-\frac{1}{5^5}$,
\begin{equation}
\Pi_{\gamma_4}(x)= \frac{5^5\sqrt{5}}{2\pi i} \left( 
t-\frac{4375}{2}t^2+\frac{16015625}{3}t^3-\frac{55322265625}{4}t^4 + 
\cdots \right),
\mylabel{eqn:Conif-t}
\end{equation}
from which $\Pi_{\gamma_4}(\frac{1}{5^5})=0$ follows (see also 
Proposition \ref{thm:gk-pi-quintic} below). 
%
The cycle $\gamma_1$ may also be understood as a `vanishing cycle' which 
vanishes at a critical point in $W^{-1}(0)$, although in this case the critical 
points are not isolated in $W^{-1}(0)$. Although the two cycles 
$\gamma_1$ and $\gamma_4$ are known explicitly, we lack such explicit
description of the cycles $\gamma_2$ and $\gamma_3$. The above Proposition 
relies only on the 
requirement that the $\gamma_k$ make a symplectic basis.
The ambiguity surrounding the description of the cycles $\gamma_2$ and $\gamma_3$
is the reason for the half integer ${\tt a}$ is left undetermined. 
However, from the following Laplace transforms, it is clear that both 
cycles $\gamma_2$ and $\gamma_3$ arise from the critical 
value $x_c=0$ of the LG potential $W(y)$ in (\ref{eqn:LGpot-W}).

\vskip0.5cm
\begin{proposition} \mylabel{thm:gk-pi-quintic} 
Applying Laplace transforms to the series (\ref{eqn:periods-quintic}), 
we obtain the formal series (\ref{eqn:formal-g-quintic}):
\begin{equation}
\begin{aligned}
&
g_k(z)=\int_0^\infty e^{-zx} \Pi_{\gamma_k}(x) dx \;\;\;\;
(k=1,2,3),  \\
&
g_4(z)=\int_{\frac{1}{5^5}}^\infty e^{-zx}\Pi_{\gamma_4}(x)dx 
=e^{-\frac{z}{5^5}} \int_0^\infty e^{-zt} \Pi_{\gamma_4}(t) dt \;\;\;. \\
\end{aligned}
\mylabel{eqn:gk-laplace}
\end{equation}
\end{proposition}

\vskip0.5cm
Now let us arrange the period integrals into a row vector 
$$
\Pi(x)=(\,\Pi_{\gamma_1}(x)\,\Pi_{\gamma_2}(x)\,\Pi_{\gamma_3}(x)\,
\Pi_{\gamma_4}(x)\,) \;\;.
$$ 
Then the monodromies about $x=0$ and $x=\frac{1}{5^5}$ are represented by 
the respective matrices acting on this row vector from the right:
$$
M_{0}=\left( \begin{matrix} 
1 & 1 & \frac{5}{2}+{\tt a} & -5 \\
0 & 1 & 5 & -\frac{5}{2}+{\tt a} \\
0 & 0 & 1 & -1 \\
0 & 0 & 0 & 1 \\ \end{matrix} \right) \;\;,\;\;M_{Con.}=
\left( \begin{matrix} 
1 & 0 & 0 & 0 \\
0 & 1 & 0 & 0 \\
0 & 0 & 1 & 0 \\
1 & 0 & 0 & 1 \\ \end{matrix} \right) 
$$
(see \cite{CdOGP} for detailed derivations of these results).
%
%
%
%
%
%
%
%
The following results are now immediate:

\vskip0.5cm
\begin{proposition} The Stokes matrices $K_R, K_{\overline R}$ in 
(\ref{eqn:connectKL-quintic}) satisfy 
\begin{equation}
K_R=M_{Con.} \;\;,\;\; K_{\overline R}\,L = M_{0} \;\;.
\mylabel{eqn:Stokes-Monod-quintic}
\end{equation}
\end{proposition}

\vskip0.5cm
The geometric interpretation of these Stokes matrices $K_{R}$ and 
$K_{\overline R}$ is closely related to the (SYZ) geometry \cite{SYZ} 
which appears at the maximal unipotent monodromy point (maximal degeneracy), 
and will be discussed in the final section. 
Here we only remark that the monodromy about $x=\infty$ is represented 
by $M_{\infty}=M_{Con.}M_{0}$, and satisfies 
$M_{\infty}^5=1$. Corresponding to these monodromy relations, we have
\begin{equation}
\begin{aligned}
&\Phi_{l+10\pi}(z) \\ 
&=\Phi_{l}(e^{-10\pi i}z) \;
K_{\overline R}\,LK_{R}L^{-1} \; 
(LK_{\overline R}\,LK_{R}L^{-2}) \;  \cdots
(L^4 K_{\overline R}\,LK_{R}L^{-5})  \\
&=\Phi_l(e^{-10\pi i}z) 
           \;K_{\overline R}\,LK_{R}\cdot K_{\overline R}\,LK_{R} \cdots
             K_{\overline R}\,LK_{R}\, L^{-5} \\
&= \Phi_l(e^{-10\pi i} z) \; L^{-5} \;\;,
\mylabel{eqn:monod-L-5}
\end{aligned}
\end{equation}
where we use repeatedly the identity 
$$
\Phi_{l+(2k+2)\pi}(z)=
\Phi_{l+2k\pi}(e^{-2\pi i} z)\;L^kK_{\overline R} L K_R L^{-k-1} \;\;,
$$
which follows from the definitions. Here, we observe a nice 
correspondence of each term in the middle line of (\ref{eqn:monod-L-5}) 
to each movement of the half line $l_-$ passing through the Stokes 
line and also the cut of the logarithms (see Fig.~2).

\vskip0.5cm

Our analysis presented above extends 
in a straightforward way to the other Calabi-Yau hypersurfaces 
studied in \cite{KT}. 
These are given by degree $d$ hypersurfaces $X(d)$ in 
$\mathbf P^4(\vec\omega)=\mathbf P^4(\omega_1,\cdots,\omega_5)$ with 
their defining data 
$$(d;\vec\omega)=(5;1^5), (6;2,1^4), (8;4,1^4), (10;5,2,1^3) \;\;.$$ 
The mirror Calabi-Yau manifold $Y(d)$ of $X(d)$ is 
again a hypersurface in a toric variety and the Hodge number 
$h^{2,1}(Y(d))=1$. In Table 1 below, we list the Stokes matrix $K_{\bar R}$ 
and also the matrices $L, N$ for each $Y(d)$. The Stokes matrix $K_{R}$ 
has a common form for all $Y(d)$ (see Proposition \ref{thm:mainI}). 
In Appendix A, we present the differential 
equations, the holomorphic solutions $\Phi_l$, and also the 
asymptotic solutions $g_k(z)$ for $d=6,8,10$.

\vskip2cm
\begin{tabular}{|c|c|c|c|} 
\hline
 $d$ & $K_{\overline R}$ & $L$ & $N$ \\
\hline
5 
& $\displaystyle{ \left( \begin{matrix} 
1 & 0 & 0 & -5 \\
0 & 1 & 0 & {\tt a}-\frac{5}{2} \\
0 & 0 & 1 & -1 \\
0 & 0 & 0 & 1 \\ \end{matrix} \right) }$
& $\displaystyle{ 
\left( \begin{matrix} 
1 & 1 & {\tt a}+\frac{5}{2} & 0 \\
0 & 1 & 5 & 0 \\
0 & 0 & 1 & 0 \\
0 & 0 & 0 & 1 \\ \end{matrix} \right) }$
& $\displaystyle{ 
\left( \begin{matrix} 
1 & -\frac{2}{5} & -\frac{2{\tt a}}{5}-2   & 1 \\
0 & \frac{2}{5} & \frac{2{\tt a}}{5}-2     &  -1 \\
0 & \frac{1}{5} & \frac{{\tt a}}{5}-\frac{1}{2} & 0 \\
0 & -\frac{1}{5} &-\frac{{\tt a}}{5}-\frac{1}{2} & 0 \\ \end{matrix} \right)} $  \\
\hline
6 
& $\displaystyle{\left(\begin{matrix} 
1 & 0 & 0 & -4 \\
0 & 1 & 0 & {\tt a}-\frac{3}{2} \\
0 & 0 & 1 & -1 \\
0 & 0 & 0 & 1 \\ \end{matrix} \right)}$ 
& $\displaystyle{\left(\begin{matrix} 
1 & 1 & {\tt a}+\frac{3}{2} & 0 \\
0 & 1 & 3 &  0 \\
0 & 0 & 1 &  0 \\
0 & 0 & 0 & 1 \\ \end{matrix} \right)}$ 
& $\displaystyle{\left(\begin{matrix} 
1 & -\frac{1}{3} & -\frac{{\tt a}}{3}-\frac{3}{2} & 1  \\
0 & \frac{1}{3}  &  \frac{{\tt a}}{3}-\frac{3}{2} &-1  \\
0 & \frac{1}{3}  &  \frac{{\tt a}}{3}-\frac{1}{2} & 0  \\
0 & -\frac{1}{3} & -\frac{{\tt a}}{3}-\frac{1}{2} & 0  \\
\end{matrix} \right)}$  \\
\hline
8 
& $\displaystyle{\left(\begin{matrix} 
1 & 0 & 0 & -4 \\
0 & 1 & 0 & {\tt a}-1 \\
0 & 0 & 1 & -1 \\
0 & 0 & 0 & 1 \\ \end{matrix} \right)}$
& $\displaystyle{\left(\begin{matrix} 
1 & 1 & {\tt a}+1 & 0 \\
0 & 1 & 2 & 0 \\
0 & 0 & 1 & 0 \\
0 & 0 & 0 & 1 \\ \end{matrix} \right)}$
& $\displaystyle{\left(\begin{matrix} 
1 & -\frac{1}{2} & -\frac{{\tt a}+3}{2} &  1 \\ 
0 & \frac{1}{2}  &  \frac{{\tt a}-3}{2} &  -1 \\ 
0 & \frac{1}{2}  &  \frac{{\tt a}-1}{2} &  0 \\ 
0 & -\frac{1}{2} & -\frac{{\tt a}+1}{2} &  0 \\ 
\end{matrix} \right)}$ \\
\hline
10 
& $\displaystyle{\left(\begin{matrix} 
1 & 0 & 0 & -3 \\
0 & 1 & 0 & {\tt a}-\frac{1}{2} \\
0 & 0 & 1 & -1 \\
0 & 0 & 0 & 1 \\ 
\end{matrix} \right)}$
& $\displaystyle{\left(\begin{matrix} 
1 & 1 & {\tt a}+\frac{1}{2} & 0 \\
0 & 1 & 1 &  0 \\
0 & 0 & 1 &  0 \\
0 & 0 & 0 & 1 \\ 
\end{matrix} \right)}$
& $\displaystyle{\left(\begin{matrix} 
1 & 0  & -1 &   1  \\
0 & 0  & -1 &  -1  \\
0 & 1  & {\tt a}-\frac{1}{2}  &  0 \\
0 & -1 & -{\tt a}-\frac{1}{2} &  0 \\
\end{matrix} \right)}$ \\
\hline
\end{tabular}
\vskip0.5cm
{\leftskip1cm \rightskip1cm 
\noindent 
{{\bf Table 1.} 
Stokes matrices $K_{\overline R}$, the monodromy matrix $L$ of the 
formal solutions $g_k(z)$ and the matrix $N$ in (\ref{eqn:fk-phi-quintic}) 
for $Y(d)$ $(d=5,6,8,10)$. The parameter ${\tt a}$ in each $d$ comes 
from the definition $g_3(z)$, and will be discussed in the final section. 
 \par} }
\vskip0.5cm

\vfill\eject 
\vskip1cm
\section{{\bf $\widehat{\bf GKZ}^*$ vs. $\widehat{\bf GKZ}_\nu$ systems}}

In Batyrev's formulation of mirror symmetry for hypersurfaces in toric 
varieties 
\cite{Ba1}, Picard-Fuchs differential 
equations  follow from the more general framework of 
Gel'fand-Kapranov-Zelevinski (GKZ) hypergeometric systems 
\cite{GKZ1,Ba2,HKTY}. 
In particular, in \cite{HKTY} it was found that GKZ systems in this context 
are resonant and reducible (and also that we need to consider a certain 
extension of the GKZ system in general). Here we consider two 
different types of Laplace transforms of the GKZ system under 
which the reducibility behaves qute differently. 

\vskip0.5cm
{\bf (3-1) Toric mirror symmetry and the GKZ system.} Let $N\simeq{\mathbf Z}^4$ 
be a lattice and $M\simeq{\mathbf Z}^4$ be its dual. An integral polytope 
in $N\otimes \mathbf R$ is called {\it reflexive} if it contains 
the origin as the only interior integral point and all its facets have 
unit integral distances. 
Given a reflexive polytope $\Delta$ in $N\otimes{\mathbf R}$, 
then its polar dual $\Delta^*$ in $M\otimes{\mathbf R}$ becomes a reflexive 
polytope, and gives rise a pair of reflexive polytopes $(\Delta,\Delta^*)$. 
In the toric mirror duality due to Batyrev, a pair $(\Delta,\Delta^*)$ 
gives a mirror pair of Calabi-Yau hypersurfaces $(X_\Delta,X_{\Delta^*})$ 
in suitable ambient toric varieties \cite{Ba1}. 
More precisely, the data of a Calabi-Yau hypersurface $X=X_\Delta$ gives 
rise to a family of mirror Calabi-Yau hypersurfaces 
$Y_{\vec a}=Y_{\Delta^*}(\vec a)$ with the defining equation,
\begin{equation}
f_{\Delta^*}(y,\vec a)=\sum_{\nu \in \Delta^*\cap M} a_\nu y^\nu \;\;,
\mylabel{eqn:LG-Q}
\end{equation}
where $\vec a \in ({\mathbf C^*})^{\Delta^*\cap M}$ denotes the (polynomial) deformations.  
The period integral we are interested in then takes the form 
\begin{equation*}
{\varPi}(\vec a)=
\frac{1}{(2\pi i)^4}
\int_\gamma Res_{f=0}\left( \frac{1}{f_{\Delta^*}(y,\vec a)}   
\prod_{k=1}^4 \frac{dy_k}{y_k} \right) \;\;,
\end{equation*}
with a choice of cycle $\gamma \in H_3(Y_{\vec a_0},\mathbf Z)$. 
Our expression $\Pi_\gamma(x)$ in (\ref{eqn:Pi[x]}) for the mirror 
quintic follows from this general expression by multiplying with a suitable 
factor. 

Let $A=\{\, (1,\nu) \,|\, \nu \in \Delta^*\cap M \,\}$ and order the elements 
of $A$ as $\bar\nu_0, \bar\nu_1, \cdots, \bar\nu_p$ $(p=|\Delta^*\cap M|-1)$ 
with $\bar\nu_k=(1,\nu_k)$ and the convention $\nu_0=(1,0)$ for 
$0 \in \Delta^* \cap M$. Then $\varPi(a)$ 
satisfies the following $A$-hypergeometric series with exponent $\beta$ 
\begin{equation*}
\square_l \varPi(\vec a)=0 \;\;(l \in L),\;\; Z \varPi(\vec a)=0 \;\;,
\end{equation*}
with
\begin{equation*}
\square_l = 
\prod_{l_i>0}\left(\frac{\pd \; }{\pd a_i} \right)^{l_i} -
\prod_{l_i<0}\left(\frac{\pd \; }{\pd a_i} \right)^{-l_i} \;\;,\;\;
Z=\sum_i \bar\nu_i \theta_{a_i} - \beta \;\;, 
\end{equation*}
where $a_i=a_{\nu_i}$ and $L\subset \mathbf Z^{p+1}$ is the lattice 
representing integral relations among the ordered elements in $A$. $Z$ 
is a differential operator taking its value in $\mathbf C^{5}$ and 
$\beta=(-1,0,0,0,0)$. 
See \cite{Ba2} for more detailed definitions. We call 
$A$-hypergoemetric systems with the exponent $\beta$ simply GKZ systems. 
The GKZ system in mirror symmetry is {\it resonant} and also {\it reducible} 
in general \cite{HKTY} (see also \cite{St}). 

In the following two sections, we consider 
two different types of Laplace transforms of the GKZ system.

\vskip0.3cm
{\bf (3-2) $\widehat{\bf GKZ}^*$ system.} The set $A$ contains a distinguished 
point $(1,0)=(1,\nu_0)$ which corresponds to the origin in $\Delta^*$. 
With respoect to the variable $a_0=a_{\nu_0}$, we define the following 
(formal) Laplace transform;
\begin{equation}
\hat\varPi(z,\bar a) = \int_0^\infty e^{-za_0} \varPi(\vec a) da_0\;\;,
\mylabel{eqn:hatPi-Q}
\end{equation}
where $\bar a=(a_1,\cdots,a_p)$. We define the  
$\widehat{\rm GKZ}^*$ system then as the Laplace transform of the GKZ system, 
obtained by the following replacements of operators in the GKZ system:
\begin{equation}
\frac{\partial \;}{\partial a_0} \rightarrow z\;,\; 
a_0 \rightarrow -\frac{\partial \;}{\partial z} \;,\; 
\theta_{a_0} \rightarrow -\theta_z-1 \;.
\mylabel{eqn:replace-a0}
\end{equation}
We note that this system is defined canonically for a given set of 
toric data of $\Delta^*$, since the origin is the unique integral point 
inside the polytope $\Delta^*$. We note further that, in toric geometry, 
the normal 
cone at the origin corresponds to the dense orbit $(\mathbf C^*)^4$ in 
the toric variety. 
It is straightforward to see that the Laplace transform 
(\ref{eqn:hatPi-Q}) represents the oscillatory integral of the Landau-Ginzburg 
theory $(f_{\Delta^*}(y,\vec a),(\mathbf C^*)^4)$ with the Landau-Ginzburg 
potential $f_{\Delta^*}(y,\vec a)$ (\ref{eqn:LG-Q}), which was first 
introduced in \cite{CV1} and studied further, for example, in 
\cite{CV2,Du1,HIV}  
in exploring the  mirror symmetry of Fano varieties.  

\vskip0.3cm
\noindent
{\bf Example 1.} 
As an illustration of the  $\widehat{\rm GKZ}^*$ system, let us 
consider the toric data $\Delta^*$ with its integral points given by 
$\nu_1=(1,0,0,0)$, 
$\nu_2=(0,1,0,0)$, 
$\nu_3=(0,0,1,0)$, 
$\nu_4=(0,0,0,1)$, 
$\nu_5=(-1,-1,-1,-1)$ 
and the origin $\nu_0=(0,0,0,0)$. This is the data for the mirror quintic $Y$
and, in fact, the Landau-Ginzburg potential coincides with (\ref{eqn:LG-f}).
The set $A$ consists of six points $(1,\nu_k)$ in $\mathbf Z^5$ and has one 
linear relation which generates the lattice $L$ of rank one. 
Under the Laplace transform, the operator $\square_l$ corresponding 
to the generator $l$ of $L$ will be replaced by 
$\hat \square_l$:
\begin{equation*}
\square_l=
\frac{\partial \;}{\partial a_1}
\frac{\partial \;}{\partial a_2}
\frac{\partial \;}{\partial a_3}
\frac{\partial \;}{\partial a_4}
\frac{\partial \;}{\partial a_5} - 
\left( \frac{\partial \;}{\partial a_0} \right)^5 
\;\mapsto \;
\hat\square_l=
\frac{\partial \;}{\partial a_1}
\frac{\partial \;}{\partial a_2}
\frac{\partial \;}{\partial a_3}
\frac{\partial \;}{\partial a_4}
\frac{\partial \;}{\partial a_5} - z^5 \;\;.
\end{equation*}
Similarly, the components of the differential operators $Z$ give rise the 
following operators
\begin{equation}
-\theta_z+\theta_{a_1}+\cdots+\theta_{a_5}\;\;,\;\;
\theta_{a_1}-\theta_{a_5}\;,\;
\theta_{a_2}-\theta_{a_5}\;,\;
\theta_{a_3}-\theta_{a_5}\;,\;
\theta_{a_4}-\theta_{a_5}\;,\;
\mylabel{eqn:hatZ-Q(5)}
\end{equation}
which annihilate $\hat\varPi(z,\bar a)$. The presence of 
the latter differential 
operators implies that $\hat\varPi(z,\bar a)$ is in fact a function of 
$\tilde z=(a_1a_2a_3a_4a_5)^{\frac{1}{5}} z$, and from $\hat \square_l$ 
we finally obtain 
\begin{equation}
\left\{ \theta_{\tilde z}^5- 5^5 \tilde z^5 \right\} 
\hat\varPi(\tilde z) = 0 \;\;.
\mylabel{eqn:Qdiff-IP4}
\end{equation}
This is the well-studied differential equation of \cite{Du2,Gu}, 
where it arose in describing the Frobenius structure of 
the quantum cohomology of 
the projective space $\mathbf P_{\Delta}=\mathbf P^4$ with 
dual polytope $\Delta$. 
\hfill \qed

\vskip0.3cm
As we see in the above example, and also observe in general, unlike
the original GKZ system, the 
$\widehat{\rm GKZ}^*$ system is {\it not} reducible. This should 
be contrasted with the second type of Laplace transform below which 
reproduces the differential equation (\ref{eqn:Qdiff-Y5}) after a 
factorization of the differential operator.

\vskip0.5cm
{\bf (3-3) $\widehat{\bf GKZ}_\nu$ system.} Now let us select one variable 
$a_\nu$ for $\nu\not=0$, and consider the second type of Laplace 
transform 
\begin{equation}
\varPi(z,\bar a)=\int_0^\infty e^{-za_\nu} \varPi(\vec a) da_\nu \;\;,
\mylabel{eqn:hatPi-Y}
\end{equation}
where $\bar a$ represents the rest of the variables in $\vec a$. 
The corresponding Laplace transform of the GKZ system is obtained by 
similar replacements to (\ref{eqn:replace-a0}) involving 
the variable $a_\nu$.  
We call the resulting set of differential equations the
$\widehat{\rm GKZ}_\nu$ system. By contrast with 
the $\widehat{\rm GKZ}^*$ system, we remark that the $\widehat{\rm GKZ}_\nu$ 
system depends on the choice of the point $\nu$ up to $GL(4,\mathbf Z)$ 
transformations on the polytope $\Delta^*$. 

\vskip0.3cm
\noindent
{\bf Example 2.} Let us consider the polytope $\Delta^*$ of the mirror quintic 
introduced in Example 1, and describe the $\widehat{\rm GKZ}_{\nu}$ 
system coming from the choice of the vertex $\nu=\nu_5$. 
In this case the choice of vertex  
is unique up to $GL(4,\mathbf Z)$ transformations. 
This time the operator $\square_l$ 
for the generator $l$ of $L$ is replaced by the following $\hat \square_l$: 
\begin{equation*}
\square_l \;\; \mapsto \;\;
\hat\square_l=
\frac{\partial \;}{\partial a_1}\frac{\partial \;}{\partial a_2}
\frac{\partial \;}{\partial a_3}\frac{\partial \;}{\partial a_4} z
- \left( \frac{\partial \;}{\partial a_0} \right)^5 \;\;.
\end{equation*}
For the components of the operator $Z$, we have 
\begin{equation}
-\theta_z+\theta_{a_0}+\cdots+\theta_{a_4}\;\;,\;\;
\theta_{a_1}+\theta_{z}+1\;,\;
\theta_{a_2}+\theta_{z}+1\;,\;
\theta_{a_3}+\theta_{z}+1\;,\;
\theta_{a_4}+\theta_{z}+1\;,\;
\mylabel{eqn:hatZ-Y(5)}
\end{equation}
which annihilate $\hat\varPi(z,\bar a)$. These relations imply that 
$\frac{a_1a_2a_3a_4}{a_0^4}\hat\varPi(z,\bar a)$ is a function of 
$\tilde z=-\frac{a_0^5}{a_1a_2a_3a_4} z$ (where the minus sign is 
simply a convention).  After some manipulations, we obtain the 
following differential equation from $\hat\square_l$:  
\begin{equation*}
\left\{ \tilde z \theta_{a_1}\theta_{a_2}\theta_{a_3}\theta_{a_4} + 
(\theta_{a_0}-4)
(\theta_{a_0}-3)
(\theta_{a_0}-2)
(\theta_{a_0}-1)
\theta_{a_0} \right\} \hat\varPi(z,\bar a)=0.
\end{equation*}
Now define 
\begin{equation*}
\hat\Pi(\tilde z)=-\frac{a_1a_2a_3a_4}{a_0^4} \hat\varPi(z, \bar a) \;\;,
\end{equation*}
then the differential equation becomes 
\begin{equation}
\theta_{\tilde z} \left\{ \tilde z (\theta_{\tilde z}+1)^3 + 
5 (5\theta_{\tilde z}+1)(5\theta_{\tilde z}+2)
(5\theta_{\tilde z}+3)(5\theta_{\tilde z}+4) \right\} \hat\Pi(\tilde z) =0 \;.
\mylabel{eqn:diffQ-Y5-GKZ}
\end{equation}
Thus we arrive at a reducible differential equation whose irreducible part 
coincides with our Laplace transform (\ref{eqn:Qdiff-Y5}). It is an easy 
exercise to derive, from (\ref{eqn:hatPi-Y}), the following oscillatory 
integral representation of $\hat\Pi(\tilde z)$, 
\begin{equation*}
\hat\Pi(\tilde z) = \frac{1}{(2\pi i)^4}\int_\Gamma e^{-\tilde z W(y)} 
dy_1dy_2dy_3dy_4 \;\;,
\end{equation*}
with the same Laudau-Ginzburg potential as (\ref{eqn:LGpot-W}). 
For this form of the oscillatory integral, one may infer that the coordinate  
$(y_1,y_2,y_3,y_4)$ is that of an affine variety. In fact, 
$y_1,\cdots,y_4$ may be identified with the coordinate 
of an affine chart of the non-compact toric variety which is defined 
by the normal cone of the vertex $\nu=\nu_5$. In other 
words, one may say that (the irreducible part of) the $\widehat{\rm GKZ}_\nu$ 
system determines the oscillatory integral of the Landau-Ginzburg theory 
$(W(y),\mathbf C^4_\nu)$ on an affine variety $\mathbf C^4_\nu$ 
determined from the vertex $\nu$. This affine Landau-Ginzburg theory 
is the one relevant to our (mirror) quintic hypersurface. 
\hfill $\square$

\vskip0.5cm
In Appendix A, we list the differential equations for the other cases of 
mirror Calabi-Yau hypersurfaces $Y(d)$. These differential equations  
also follow from the $\widehat{\rm GKZ}_\nu$ systems under a suitable 
choice of the vertex $\nu$.

\vskip1cm
\section{{\bf Conclusion and discussions}}
\vskip0.5cm

By making Laplace transforms of Picard-Fuchs equations of period integrals, 
we have obtained differential equations with irregular singularities 
at infinity which are degenerate. We have determined the Stokes matrices 
of the differential equations. The same differential equations have been 
obtained by introducing the $\widehat{\rm GKZ}_\nu$ system, a Laplace transform 
of GKZ system. We contrasted the $\widehat{\rm GKZ}_\nu$ system with a similar 
Laplace transform $\widehat{\rm GKZ}^*$, which reproduces the differential 
equations related to the quantum cohomology of toric varieties. 

Although in this note we have restricted our attention to the simplest 
cases of Calabi-Yau hypersurfaces with one parameter deformations, similar 
analyses should be possible for more general Calabi-Yau hypersurfaces starting 
with our $\widehat{\rm GKZ}_\nu$ system. 
Also similar Laplace transforms may be considered 
for Calabi-Yau varieties such as complete intersections 
in Grassmanians, and in particular for the 4th order differential operators 
of Calabi-Yau type classified in \cite{DM,ES}. However, 
in such generalizations, Landau-Ginzburg theories become less clear 
than for the hypersurface case.  

\vskip0.5cm
In the rest of this section, we present an interpretation of our Stokes 
matrices from the perspective of mirror symmetry between $X(d)$ and $Y(d)$. 

\vskip0.5cm

In the papers \cite{CdOGP,KT}, the parameter ${\tt a}$ in the Stokes 
matrices has been set to ${\tt a}=\frac{11}{2}, \frac{9}{2}, 3, \frac{1}{2}$ 
for $d=5,6,8,10$, respectively. This parameter is directly related to the 
ambiguity in the quadratic term of the prepotential (the gererating 
function of the Gromov-Witten invariants) and hence plays no role in 
the study of Gromov-Witten invariants. 
From our viewpoint of vanishing cycles, however, we should be able to 
fix this parameter by defining and integrating periods over 
the symplectic cycles $\gamma_2,\gamma_3$. 
Here, from Table 1, we simply observe that by choosing 
${\tt a}=\frac{5}{2}, \frac{3}{2}, 1, \frac{1}{2}$ for $d=5,6,8,10$, 
respectively, the Stokes matrix $K_{\overline R}$ simplifies. Furthermore, 
we may identify the form of the Stokes matrices in terms of 
the invariants of $X(d)$ as follows 
\begin{equation*}
K_R=\begin{pmatrix} 
1 & 0 & 0 &  0 \\
0 & 1 & 0 &  0  \\
0 & 0 & 1 &  0  \\
\chi(\mathcal O_X,\mathcal O_p) & 0 & 0 &  1  \\
\end{pmatrix} \;,\;\;
K_{\overline R}=\begin{pmatrix} 
1 & 0 & 0 & -\chi(\mathcal O_X, \mathcal O_X(1) ) \\
0 & 1 & 0 &  0  \\
0 & 0 & 1 & -1  \\
0 & 0 & 0 &  1  \\
\end{pmatrix} \;,
\mylabel{eqn:Stokes-chi}
\end{equation*}
where $\mathcal O_X$ and $\mathcal O_p$ represent the structure sheaf of 
$X$ and a skyscraper sheaf on a point $p$ of $X$, respectively, and 
$\mathcal O_X(1)=\mathcal O_X \otimes H$, with $H$ the ample line bundle 
which generates the Picard group ${\rm Pic}(X)=\mathbf Z H$. 
Also, for coherent 
sheaves $\mathcal E, \mathcal F$ on $X$
$$\chi(\mathcal E,\mathcal F)=\sum_{k=0}^3 (-1)^k 
\dim {\rm Hom}(\mathcal E,\mathcal F) \;\;.$$ 
By the Riemann-Roch theorem, one may easily verify that  
$\chi(\mathcal O_X, \mathcal O_X(1)) = 5, 4, 4, 3,$
for $X(d)\; (d=5, 6, 8, 10)$, respectively. The fact that the entry $-1$ in 
$K_{\bar R}$ is common for all $Y(d)$ suggests a similar invariant 
interpretation, but we shall defer this to future studies. 
Clearly, these invariant expressions are explained 
by the homological mirror symmetry conjecture, 
which implies the correspondence of  
$\mathcal O_X$ and $\mathcal O_p$ with the vanishing cycle $\gamma_4=S^3$ and 
the torus cycle $\gamma_1=T^3$, respectively \cite{SYZ,Gr}. 
Then, $K_{R}$ represents the 
Picard-Lefschetz transform for the vanishing cycle $\gamma_4$, while 
$K_{\overline R}$ corresponds to the case where the 
torus cycle $\gamma_1$ vanishes. We can see these degenerations at the 
critical loci in $W^{-1}(\frac{1}{5^5})$ and $W^{-1}(0)$, respectively, 
of the the Landau-Ginzburg potential 
\begin{equation*}
W(y)=y_1y_2y_3y_4(y_1+y_2+y_3+y_4+1) \;\;.
\end{equation*}
In $W^{-1}(\frac{1}{5^5})$, we see that the cycle $\gamma_4\approx S^4$ 
vanishes 
and observe the Picard-Lefshetz transformation in the Stokes matrix $K_R$. 
On the other hand, $W^{-1}(0)$ has a degeneration charaterized by  
maximal unipotent monodromy.  Near this degeneration, 
a torus fibration, called the SYZ fibration \cite{SYZ}, has been 
constructed in 
\cite{Zh,Ru} using the amoeba picture \cite{GKZ2,Mi} of 
a suitable moment map image of $W^{-1}(\epsilon)$. When $\epsilon$ goes to 
zero, the moment map image $W(y)=\epsilon$ converges to the boundary of a 
four dimensional simplex which roughly describes the cycle 
$\gamma_4 \approx S^3$. The dual cycle $\gamma_1
\approx T^3$ vanishes at each point on the (five) faces of the 
simplex, which explains the number $-5$ in $K_{\overline R}$. In this picture, 
the `Picard-Lefshetz' transform $K_{\overline R}$ represents a twisting of 
the zero section $\gamma_4$ (cf. \cite{Ab}). 
Our missing cycles $\gamma_2,\gamma_3$ in the symplectic basis 
should be found in the geometry of $W^{-1}(0)$ in a similar way.

It is worthwhile reproducing here the well-known Stokes matrices 
$S, S_-$ of the differential equation (\ref{eqn:Qdiff-IP4}) 
for the quantum cohomology 
of the projective space $\mathbf P^4$ \cite{Du2,Gu}. Under a suitable 
choice of basis, the Stokes matrices are calculated in \cite{Gu} 
and found to be $S=
(\chi(\mathcal O_{\mathbf P^4}(i),
 O_{\mathbf P^4}(j)) )_{0\leq i,j \leq 4}^{-1}$ 
and $S_-=S^T$ in terms of the coherent sheaves 
$\mathcal O_{\mathbf P^4}, \mathcal O_{\mathbf P^4}(1), 
\cdots,\mathcal O_{\mathbf P^4}(4)$ on $\mathbf P^4$. In this case, 
the Landau-Ginzburg potential $f(y,a)$ in (\ref{eqn:LG-f}) 
has five isolated critical points of ordinary double points, and each is 
considered to be the mirror of $\mathcal O_{\mathbf P^4}(i)$ under a suitable 
ordering. Thus, we see that although the $\widehat{\rm GKZ}^*$ and 
$\widehat{\rm GKZ}_\nu$ systems are defined in a similar way starting from the 
same GKZ system, they have completely different geometric interpretations 
in mirror symmetry. 

\vskip0.5cm

\vfill\eject
\appendix
\section{{\bf Calabi-Yau hypersurfaces in $\mathbf P^4(\omega)$}}

According to \cite{KT}, we define Calabi-Yau hypersurfaces $X(d)$ in 
weighted projective spaces $\mathbf P^4(\vec\omega)=
\mathbf P^4(\omega_1,\cdots,\omega_5)$ for each set of toric data 
$$(d;\vec\omega) = (5;1^5), (6;2,1^4), (8;4,1^4), (10;5,2,1^3) \;\;.$$ 
The mirror $Y(d)$ of $X(d)$ has Hodge number $h^{2,1}(Y(d))=1$. 
The mirror pairs $(X(d),Y(d))$ are typical examples of Batyrev's mirror 
symmetry associated to the corresponding pairs of reflexive polytopes 
$(\Delta(\vec\omega),\Delta(\vec\omega)^*)$. The Laplace transforms 
of the Picard-Fuchs differential equations listed below follow from 
the $\widehat{\rm GKZ}_\nu$ under a suitable choice of a vertex $\nu$.

\vskip0.3cm
{\bf (A-1) Series data and Stokes matrices of $Y(d)$.} 
$Y(5)$ is the quntic mirror hypersurface which has been studied in the 
text. Since analyses of other cases $Y(d)$ are parallel, we 
list the corresponding results below, where the items (a) and (b) 
are the series data near $z=0$ and $z=\infty$, respectively. For 
convenience, we have included the Picard-Fuchs equations as well as their 
Laplace transforms. 

\vskip0.3cm
\noindent
\underline{$\bullet$ $Y(6)$:}
The Picard-Fuchs equation and its Laplace transform, respectively, are 
\noindent
\begin{equation*}
\begin{aligned}
 \{ \theta_x^4-9x(6\theta_x+5)(6\theta_x+4)(6\theta_x+2)(6\theta_x+1)\}
\Pi(x)&=0\;, \\
\{ z(\theta_z+1)^3+9(6\theta_z+1)(6\theta_z+2)(6\theta_z+4)(6\theta_z+5) \}
\hat \Pi(z) &=0 \;. 
\end{aligned}
\end{equation*}

\noindent
(a) We set 
\begin{equation*}
\phi_{k}(z)=\sum_{n=0}^\infty 
\frac{\Gamma(1+n+\rho_k)^3\Gamma(3+2(n+\rho_k))}{\Gamma(7+6(n+\rho_k))}
(-1)^n z^{n+\rho_k} \;,
\quad 
\end{equation*}
where $\rho_k=-\frac{k}{6}$ for $k=1,2,4,5$. 
We may write these solutions by Barnes hypergeometric series
\begin{equation*}
\begin{aligned}
&
\phi_1(z)=
\frac{d_1}{z^{\frac{1}{6}}} 
\,_3F_3(\frac{5}{6}^3;\frac{7}{6},\frac{3}{2},\frac{5}{3},
\frac{-z}{2^4 3^6}) ,\;\;
\phi_2(z)=
\frac{d_2}{z^{\frac{2}{6}}}
\,_3F_3(\frac{2}{3}^3;\frac{5}{6},\frac{4}{3},\frac{3}{2},
\frac{-z}{2^4 3^6}) \;\;,\;\; \\
&
\phi_4(z)=
\frac{d_4}{z^{\frac{4}{6}}} 
\,_3F_3(\frac{1}{3}^3;\frac{1}{2},\frac{2}{3},\frac{7}{6},
\frac{-z}{2^4 3^6}) ,\;\;
\phi_5(z)=
\frac{d_5}{z^{\frac{5}{6}}}
\,_3F_3(\frac{1}{6}^3;\frac{1}{3},\frac{1}{2},\frac{5}{6},
\frac{-z}{2^4 3^6}) \;\;,\;\; \\
\end{aligned}
\end{equation*}
where the constants $d_1,d_2,d_4,d_5$ are given by 
$\frac{\Gamma(\frac{2}{3})\Gamma(\frac{5}{6})}{2^2 3^3}$, 
$\frac{\Gamma(\frac{1}{3})\Gamma(\frac{2}{3})^3}{54}$, 
$\frac{\Gamma(\frac{1}{3})^3\Gamma(\frac{5}{3})}{2}$, and 
$\Gamma(\frac{1}{6})^3\Gamma(\frac{4}{3})$, respectively.
Define 
\begin{equation*}
\hat w^{\infty}_m(z)=  (1-\alpha)
\sum_{k=1}^5 (1-\alpha^{2k})(1-\alpha^k)^2 
\alpha^{km+\frac{1}{2}(k-1)}\,\phi_{6-k}(z) ,
\end{equation*}
for $0\leq m \leq 5$ with $\alpha=e^{2\pi i/6}$. Then the holomorphic 
solution is given by 
\begin{equation*}
\Phi_l(z)=c_N( \;
\hat w^{\infty}_0(z) \; \hat w^{\infty}_1(z) \; 
\hat w^{\infty}_2(z) \; \hat w^{\infty}_5(z) \; ) N \;\;,
\end{equation*}
where $c_N=-0.00201572...$ and the matrix $N$ is given in Table 1. 

\noindent
(b) We set 
\begin{equation*}
G(z,\rho)=\sum_{n=0}^\infty 
\frac{\Gamma(1+6(n+\rho))}{\Gamma(1+n+\rho)^3\Gamma(1+2(n+\rho))} 
\left( \frac{1}{z} \right)^{n+\rho+1} \;\;,
\end{equation*}
and define the asymptotic solutions 
\begin{equation*}
\begin{aligned}
& 
g_1(z)=G(z,0) \;,\;\; 
g_2(z)=\pd_{\tilde\rho}G(z,\rho)\Big|_{\rho=0} ,\;\;
g_3(z)=\big\{ \frac{3}{2} \pd_{\tilde\rho}^2+{\tt a} \pd_{\tilde\rho} \big\} 
G(z,\rho)\Big|_{\rho=0} \\
&
g_4(z)=\frac{2^43^6\sqrt{3}}{2\pi i} e^{-\frac{z}{2^4 3^6}} \frac{1}{z^2} 
\left(1-\frac{15876}{z}+\frac{428354568}{z^2}-\cdots \right)\;, \\
\end{aligned}
\mylabel{eqn:formal-g-X(6)}
\end{equation*}
where $\pd_{\tilde\rho}=\frac{1}{2\pi i} \frac{\pd \;}{\pd \rho}$ and 
${\tt a}$ is a parameter.

\vskip0.3cm
\noindent
\underline{$\bullet$ $Y(8)$:} 
The Picard-Fuchs equation and its Laplace transform, respectively, are
\begin{equation*}
\begin{aligned}
\{ \theta_x^4-16x(8\theta_x+7)(8\theta_x+5)(8\theta_x+3)(8\theta_x+1)\}
\Pi(x)&=0\;, \\
\{ z(\theta_z+1)^3+16(8\theta_z+1)(8\theta_z+3)(8\theta_z+5)(8\theta_z+7) \}
\hat \Pi(z) &=0 \;.
\end{aligned}
\end{equation*}

\noindent
(a) We set 
\begin{equation*}
\phi_{k}(z)=\sum_{n=0}^\infty 
\frac{\Gamma(1+n+\rho_k)^3\Gamma(5+4(n+\rho_k))}{\Gamma(9+8(n+\rho_k))}
(-1)^n z^{n+\rho_k} 
\quad 
\end{equation*}
where $\rho_k=-\frac{k}{8}$ for $k=1,3,5,7$. 
We may write these solutions by Barnes hypergeometric series
\begin{equation*}
\begin{aligned}
&
\phi_1(z)=
\frac{d_1}{z^{\frac{1}{8}}} 
\,_3F_3(\frac{7}{8}^3;\frac{5}{4},\frac{6}{4},\frac{7}{4},
\frac{-z}{2^{16}}) ,\;\;
\phi_3(z)=
\frac{d_3}{z^{\frac{3}{8}}}
\,_3F_3(\frac{5}{8}^3;\frac{3}{4},\frac{5}{4},\frac{6}{4},
\frac{-z}{2^{16}}) \;\;,\;\; \\
&
\phi_5(z)=
\frac{d_5}{z^{\frac{5}{8}}} 
\,_3F_3(\frac{3}{8}^3;\frac{2}{4},\frac{3}{4},\frac{5}{4},
\frac{-z}{2^{16}}) ,\;\;
\phi_7(z)=
\frac{d_7}{z^{\frac{7}{8}}}
\,_3F_3(\frac{1}{8}^3;\frac{1}{4},\frac{2}{4},\frac{3}{4},
\frac{-z}{2^{16}}) \;\;,\;\; \\
\end{aligned}
\end{equation*}
where the constants $d_1,d_3,d_5,d_7$ are given by 
$\frac{\Gamma(\frac{1}{2})\Gamma(\frac{7}{8})^3}{2^8  3}$, 
$\frac{\Gamma(\frac{1}{2})\Gamma(\frac{5}{8})^3}{2^6}$, 
$\frac{\Gamma(\frac{1}{2})\Gamma(\frac{3}{8})^3}{2^3}$, and 
$\frac{\Gamma(\frac{1}{2})\Gamma(\frac{1}{8})^3}{2}$, respectively.
Define
\begin{equation*}
\hat w^{\infty}_m(z)=  (1-\alpha)
\sum_{k=1}^7 (1-\alpha^{4k})(1-\alpha^k)^2 
\alpha^{km+\frac{1}{2}(k-1)}\,\phi_{8-k}(z) ,\;\; 
\end{equation*}
for $0\leq m \leq 7$ with $\alpha=e^{2\pi i/8}$. Then the holomorphic 
solutions are 
\begin{equation*}
\Phi_l(z)=c_N( \;
\hat w^{\infty}_0(z) \; \hat w^{\infty}_1(z) \;
\hat w^{\infty}_2(z) \; \hat w^{\infty}_7(z) \;) N \;\;,
\end{equation*}
where $c_N=-0.001316833...$ and the matrix $N$ is given in Table 1.

\noindent
(b) We set 
\begin{equation*}
G(z,\rho)=\sum_{n=0}^\infty 
\frac{\Gamma(1+8(n+\rho))}{\Gamma(1+n+\rho)^3\Gamma(1+4(n+\rho))} 
\left( \frac{1}{z} \right)^{n+\rho+1} \;\;,
\end{equation*}
and define 
\begin{equation*}
\begin{aligned}
& 
g_1(z)=G(z,0) \;,\;\; 
g_2(z)=\pd_{\tilde\rho}G(z,\rho)\Big|_{\rho=0} ,\;\;
g_3(z)=\big\{ \frac{2}{2} \pd_{\tilde\rho}^2+{\tt a} \pd_{\tilde\rho} \big\} 
G(z,\rho)\Big|_{\rho=0} \\
&
g_4(z)=\frac{2^{16}\sqrt{2}}{2\pi i} e^{-\frac{z}{2^{16}}} \frac{1}{z^2} 
\left(1-\frac{88064}{z}+\frac{13272875008}{z^2}-\cdots \right)\;, \\
\end{aligned}
\mylabel{eqn:formal-g-X(8)}
\end{equation*}
where $\pd_{\tilde\rho}=\frac{1}{2\pi i} \frac{\pd \;}{\pd \rho}$ and 
${\tt a}$ is a parameter.

%

\vskip0.3cm

\noindent
\underline{$\bullet$ $Y(10)$:} 
The Picard-Fuchs equation and its Laplace transform, respectively, are
\begin{equation*}
\begin{aligned}
\{ \theta_x^4-2^4 5x(10\theta_x+9)(10\theta_x+7)(10\theta_x+3)(10\theta_x+1)\}
\Pi(x)&=0\;, \\
\{ z(\theta_z+1)^3+2^4 5(10\theta_z+1)(10\theta_z+3)(10\theta_z+7)(10\theta_z+9) \}
\hat \Pi(z)&=0 \;.
\end{aligned}
\end{equation*}

\noindent
(a) We set 
\begin{equation*}
\phi_{k}(z)=\sum_{n=0}^\infty 
\frac{\Gamma(1+n+\rho_k)^2\Gamma(3+2(n+\rho_k))\Gamma(6+5(n+\rho_k))}
{\Gamma(11+10(n+\rho_k))}
(-1)^n z^{n+\rho_k} 
\quad 
\end{equation*}
where $\rho_k=-\frac{k}{10}$ for $k=1,3,7,9$. 
We may write these solutions by Barnes hypergeometric series
\begin{equation*}
\begin{aligned}
&
\phi_1(z)=
\frac{d_1}{z^{\frac{1}{10}}} 
\,_3F_3(\frac{9}{10}^3;\frac{5}{6},\frac{8}{5},\frac{9}{5},
\frac{-z}{2^8 5^5}) ,\;\;
\phi_3(z)=
\frac{d_3}{z^{\frac{3}{10}}}
\,_3F_3(\frac{7}{10}^3;\frac{4}{5},\frac{7}{5},\frac{8}{5},
\frac{-z}{2^8 5^5}) \;\;,\;\; \\
&
\phi_7(z)=
\frac{d_7}{z^{\frac{7}{10}}} 
\,_3F_3(\frac{3}{10}^3;\frac{2}{5},\frac{3}{5},\frac{6}{5},
\frac{-z}{2^8 5^5}) ,\;\;
\phi_9(z)=
\frac{d_9}{z^{\frac{9}{10}}}
\,_3F_3(\frac{1}{10}^3;\frac{1}{5},\frac{2}{5},\frac{4}{5},
\frac{-z}{2^8 5^5}) \;\;,\;\; \\
\end{aligned}
\end{equation*}
where the constants $d_1,d_3,d_7,d_9$ are given by 
$\frac{\Gamma(\frac{9}{10})^2\Gamma(\frac{4}{5})\Gamma(\frac{1}{2})}
{2^{10} 15}$, 
$\frac{\Gamma(\frac{7}{10})^2\Gamma(\frac{2}{5})\Gamma(\frac{1}{2})}
{2^7 15}$, 
$\frac{\Gamma(\frac{3}{10})^2\Gamma(\frac{3}{5})\Gamma(\frac{1}{2})}
{2^3}$, and 
$\frac{\Gamma(\frac{1}{10})^2\Gamma(\frac{1}{5})\Gamma(\frac{1}{2})}
{2}$, respectively. 
Define 
\begin{equation*}
\hat w^{\infty}_m(z)=  (1-\alpha)
\sum_{k=1}^{9} (1-\alpha^{5k})(1-\alpha^{2k})(1-\alpha^k) 
\alpha^{km+\frac{1}{2}(k-1)}\,\phi_{10-k}(z) ,\;\; 
\end{equation*}
for $0\leq m \leq 9$ with $\alpha=e^{2\pi i/10}$. 
Then the holomorphic solutions are given by 
\begin{equation*}
\Phi_l(z)=c_N( \;
\hat w^{\infty}_0(z) \; \hat w^{\infty}_1(z) \;
\hat w^{\infty}_2(z) \; \hat w^{\infty}_9(z) \;) N \;\;,
\end{equation*}
where $c_N=-0.0001304601...$ and the matrix $N$ is given in Table 1. 

\noindent
(b) We set 
\begin{equation*}
G(z,\rho)=\sum_{n=0}^\infty 
\frac{\Gamma(1+10(n+\rho))}{\Gamma(1+n+\rho)^2\Gamma(1+2(n+\rho))
\Gamma(1+5(n+\rho))} \left( \frac{1}{z} \right)^{n+\rho+1} \;\;,
\end{equation*}
and define the asymptotic solutions 
\begin{equation*}
\begin{aligned}
& 
g_1(z)=G(z,0) \;,\;\; 
g_2(z)=\pd_{\tilde\rho}G(z,\rho)\Big|_{\rho=0} ,\;\;
g_3(z)=\big\{ \frac{1}{2} \pd_{\tilde\rho}^2+{\tt a} \pd_{\tilde\rho} \big\} 
G(z,\rho)\Big|_{\rho=0} \\
&
g_4(z)=\frac{2^{8}5^{5}}{2\pi i} e^{-\frac{z}{2^{8}5^5}} \frac{1}{z^2} 
\left(1-\frac{1040000}{z}+\frac{1884800000000}{z^2}-\cdots \right)\;, \\
\end{aligned}
\mylabel{eqn:formal-g-X(10)}
\end{equation*}
where $\pd_{\tilde\rho}=\frac{1}{2\pi i} \frac{\pd \;}{\pd \rho}$ and 
${\tt a}$ is a parameter. 

\vskip0.5cm



\vskip0.3cm
{\bf (A-2) Stokes matrices and monodromy matrices.} As in Proposition 
\ref{thm:gk-pi-quintic}, there are relations between the formal 
solutions $g_k(z)$ and period integrals $\Pi_{\gamma_k}(x)$ for 
a symplectic basis $\gamma_1,\cdots,\gamma_4$ of $H_3(Y(d),\mathbf Z)$, 
with $|\gamma_1\cap\gamma_4|=|\gamma_2\cap\gamma_3|=1$ and vaninishing 
intersections for other pairings.  Similar to the case $Y(5)$, we  
know the cycles $\gamma_1 \approx T^3$ and $\gamma_4 \approx S^3$ explicitly. 
However, as for $\gamma_2$ and $\gamma_3$, we do not have precise forms of 
the cycles. Here we only verify the expected 
symplectic and integral monodromy of the hypergeometric series 
$\Pi_{\gamma_k}(x)$ defined below\cite{KT,Ho2,DM}. 

To define $\Pi_{\gamma_k}(x)$ in general, let us set 
\begin{equation*}
w(x,\rho)=\sum_{n=0}^{\infty} \frac{\Gamma(1+d(n+\rho))}{
\prod_{i=1}^{5} \Gamma(1+\omega_i(n+\rho))} x^{n+\rho} \;\;,
\end{equation*}
for each $Y(d)$ with toric data $(d;\vec \omega)$. Then the symplectic 
bases given in Proposition \ref{thm:sympPi} generalize to $Y(d)$ by 
\begin{equation}
\begin{aligned}
&\Pi_{\gamma_1}(x) = w(x,0) \;\;,\;\;  
&\Pi_{\gamma_3}(x)=\big\{ \frac{K_d}{2}\pd_{\tilde\rho}^2 + 
               {\tt a} \pd_{\tilde\rho} \big\} w(x,\rho)\big|_{\rho=0} \;\;, 
   \hskip0.5cm \\
&\Pi_{\gamma_2}(x)=\pd_{\tilde\rho} w(x,\rho)\big|_{\rho=0} \;\;,\;\;
&\Pi_{\gamma_4}(x)=\big\{ - \frac{K_d}{3!}\pd_{\tilde\rho}^3 
      - \frac{C_d}{12} \pd_{\tilde\rho} \big\} w(x,\rho)\big|_{\rho=0} \;\;, 
   \\
\end{aligned}
\mylabel{eqn:periods-Y(d)}
\end{equation}
with $(K_d,C_d)=(5,50),(3,42),(2,44),(1,34)$ for $d=5,6,8,10$, respectively.

Analytic continuation of $\Pi_{\gamma_4}(x)$ about the relevant conifold point has 
a similar form to (\ref{eqn:Conif-t}). Then the suitable Laplace transforms of 
$\Pi_k(x)$, see (\ref{eqn:gk-laplace}), give the asymptotic series $g_k(z)$. 
For all cases of mirror hypersurfaces $Y(d)$, one verifies the relations 
\begin{equation*}
K_R=M_{Con.} \;\;,\;\; K_{\overline R} L = M_{0} \;\;,
\end{equation*}
which connect the monodromy of $\Pi_{\gamma_k}(x)$ to the Stokes 
matrices of the holomorphic solution $\Phi_l(z)$.

\vskip0.5cm

\vskip2cm

\end{document}